\RequirePackage{fix-cm}
\documentclass[envcountsect]{svjour3}                     

\smartqed  
\usepackage{graphicx}
\usepackage{caption}
\usepackage{amsfonts}
\usepackage{mathptmx}

\usepackage[]{geometry}
\usepackage[square,numbers,sort&compress]{natbib}
\spnewtheorem{fact}{Fact}[section]{\bf}{\it}

\begin{document}
\title{On Lagrangians of Hypergraphs Containing Dense Subgraphs
}


\author{Qingsong Tang         \and   Yuejian  Peng  \and   Xiangde Zhang \and Cheng Zhao
}


\institute{Qingsong Tang \at
              College of Sciences, Northeastern University, Shenyang, 110819, P.R.China. \at
              School of Mathematics, Jilin University, Changchun 130012, P.R. China.\\
              \email{t\_qsong@sina.com}
              \and
              Yuejian  Peng \at
             School of Mathematics, Hunan University, Changsha 410082, P.R. China. \at This research is supported by National Natural Science Foundation of China (No. 11271116).\\
              \email{ypeng1@163.com}
               \and
           Xiangde Zhang \at
             College of Sciences, Northeastern University, Shenyang, 110819, P.R.China\\
             \email{zhangxdneu@163.com}
              \and
             Cheng Zhao(Corresponding author) \at
             Department of Mathematics and Computer Science, Indiana State University, Terre Haute, IN, 47809 USA. \at School of Mathematics, Jilin University, Changchun 130012, P.R. China.\\
             \email{cheng.zhao@indstate.edu}
}

\date{Received: date / Accepted: date}

\maketitle

\begin{abstract}
Motzkin and Straus established a remarkable connection between the maximum
clique and the Lagrangian of a graph in 1965. This connection and its extensions were successfully
employed in optimization to provide heuristics for the maximum clique number in
graphs. It is useful in practice if similar results hold for hypergraphs. In this paper, we provide
upper bounds on the Lagrangian of a  hypergraph containing dense subgraphs
when the number of edges of the hypergraph is in certain ranges. These results support a pair  of conjectures
introduced by Y. Peng and C. Zhao (2012) and extend a result of J. Talbot (2002).
\keywords{Cliques of hypergraphs \and Colex ordering \and Lagrangians of hypergraphs \and Polynomial optimization}
 \subclass{05C35 \and 05C65 \and 05D99 \and 90C27}
\end{abstract}

\section{Introduction}
\label{intro}
In 1941, Tur\'an \cite{Turan} provided an answer to the following question: What is the maximum
number of edges in a graph with n vertices not containing a complete subgraph of order k, for
a given k? This is the well-known Tur\'an theorem. Later, in another classical paper, Motzkin
and Straus \cite{MS} provided a new proof of Tur\'an theorem based on the  continuous characterization
of the clique number of a graph using Lagrangians of graphs.

The Motzkin-Straus result basically says that the Lagrangian of a  graph which is the maximum of
a  homogeneous quadratic multilinear function (determined by the graph) over the standard simplex of the Euclidean  plane is connected to
the maximum clique number of this graph (the precise statement is given in Theorem \ref{MStheo}). This result provides a solution to the optimization problem for a class of homogeneous quadratic multilinear functions over the standard simplex of an Euclidean  plane.
The Motzkin-Straus result
and its extension were successfully employed in optimization to provide heuristics for the
maximum clique problem \cite{B1,B2,B3,G9}.  It has  been also generalized to vertex-weighted graphs \cite{G9} and edge-weighted graphs with applications to pattern recognition in image analysis \cite{B1,
B2, B3, G9, PP, PP15, RTP20}
The Lagrangian of a hypergraph has been a useful tool in hypergraph extremal problems.
For example, Sidorenko \cite{sidorenko89} and Frankl-Furedi \cite{FF} applied Lagrangians of hypergraphs in
finding Tur\'an densities of hypergraphs. Frankl and R\"{o}dl \cite{FR84} applied it in disproving Erd\"os
long standing jumping constant conjecture.
In most applications, we need an upper bound for the Lagrangian of a
hypergraph.

An attempt to generalize the Motzkin-Straus theorem to hypergraphs is due to S\'os and
Straus\cite{SS}. Recently, in \cite{BP1,BP2} Rota Bul\'o and Pelillo generalized the Motzkin and Straus'
result to $r$-graphs in some way using a continuous characterization of maximal cliques
other than Lagrangians of hypergraphs. The obvious generalization of Motzkin and Straus'
result to hypergraphs is false. In fact, there are many examples of hypergraphs that do not
achieve their Lagrangian on any proper subhypergraph. We attempt to explore the relationship
between the Lagrangian of a hypergraph and the order of its maximum cliques for
hypergraphs when the number of edges is in certain ranges though the obvious generalization
of Motzkin and Straus' result to hypergraphs is false.

The results presented in Sect. 3 and 4 in this paper provide substantial evidence for two conjectures in \cite{PZ}
and extend some known results in the literature \cite{PZ,T}. The main results provide solutions to the optimization problem of a class of homogeneous  multilinear functions  over the standard simplex of the  Euclidean  space. The main results also give connections between a continuous optimization problem  and the maximum clique problem of hypergraphs. Since  practical problems such as computer vision and image analysis are related to the maximum clique problems, this type of results opens a  door to  such practical applications.  The  results in this paper can be applied in estimating Lagrangians of some hypergraphs, for example, calculations involving estimating Lagrangians of several hypergraphs in \cite{FF} can be much simplified when applying the  results in this paper.

The rest of the paper is
organized as follows. In Sect. 2, we state a few definitions, problems, and preliminary
results. In Sect. 3 and Sect. 4, we provide upper bounds on the Lagrangian of a
hypergraph containing dense subgraphs when the number of edges of the hypergraph is in
a certain range. Then, as an application, using the main result in Sect. 3, we extend a result
in \cite{T} in Sect. 5. In Sect. 6, we give the proofs of some lemmas. Conclusions are given in Section 7.

\section{Definitions and Preliminary  Results}
\label{sec:1}
For a set $V$ and a positive integer $r$ we denote by $V^{(r)}$ the family of all $r$-subsets of $V$.
An $r$-uniform graph or $r$-graph $G$ consists of a set $V(G)$ of vertices and a set $E(G) \subseteq V(G) ^{(r)}$ of edges. An edge $e:=\{a_1, a_2, \ldots, a_r\}$ will be simply denoted by $a_1a_2 \ldots a_r$.
An $r$-graph $H$ is  a {\it subgraph} of an $r$-graph $G$, denoted by $H\subseteq G$ if $V(H)\subseteq V(G)$ and $E(H)\subseteq E(G)$.   Let $K^{(r)}_t$ denote the complete $r$-graph on $t$ vertices, that is the $r$-graph on $t$ vertices containing all possible edges. A complete $r$-graph on $t$ vertices is also called a clique with order $t$.
Let ${\mathbb N}$ be the set of all positive integers. For any integer $n \in {\mathbb N}$, we denote the set $\{1, 2, 3, \ldots, n\}$ by $[n]$. Let $[n]^{(r)}$  represent the  complete $r$-uniform graph on the vertex set $[n]$.
When $r=2$, an $r$-uniform graph is a simple graph.  When $r\ge 3$,  an $r$-graph is often called a hypergraph.

For an $r$-graph $G=(V,E)$ and $i\in V$, let $E_i:=\{A \in V^{(r-1)}: A \cup \{i\} \in E\}$. For a pair of vertices $i,j \in V$, let $E_{ij}:=\{B \in V^{(r-2)}: B \cup \{i,j\} \in E\}$. Let

$E^c_i:=\{A \in V^{(r-1)}: A \cup \{i\} \in V^{(r)} \backslash E\}$, $E^c_{ij}:=\{B \in V^{(r-2)}: B \cup \{i,j\} \in V^{(r)} \backslash E\}$, and $E_{i\setminus j}:=E_i\cap E^c_j.$


\begin{definition} \label{definitionlagrangian}
For  an $r$-uniform graph $G$ with the vertex set $[n]$, edge set $E(G)$, and a vector \\$\vec{x}:=(x_1,\ldots,x_n) \in {\mathbb R}^n$,
we associate a homogeneous polynomial in $n$ variables, denoted by  $\lambda (G,\vec{x})$  as follows:
$$\lambda (G,\vec{x}):=\sum_{i_1i_2 \cdots i_r \in E(G)}x_{i_1}x_{i_2}\ldots x_{i_r}.$$
Let $S:=\{\vec{x}:=(x_1,x_2,\ldots ,x_n): \sum_{i=1}^{n} x_i =1, x_i
\ge 0 {\rm \ for \ } i=1,2,\ldots , n \}$.
Let $\lambda (G)$ represent the maximum
 of the above homogeneous  multilinear polynomial of degree $r$ over the standard simplex $S$. Precisely
 $$\lambda (G): = \max \{\lambda (G, \vec{x}): \vec{x} \in S \}.$$
\end{definition}
The value $x_i$ is called the {\em weight} of the vertex $i$.
A vector $\vec{x}:=(x_1, x_2, \ldots, x_n) \in {\mathbb R}^n$ is called a feasible weighting for $G$ if
$\vec{x}\in S$. A vector $\vec{y}\in S$ is called an {\em optimal weighting} for $G$
if $\lambda (G, \vec{y})=\lambda(G)$.

\begin{remark}
Since $\lambda (G)$ is the maximum of a   polynomial function  in $n$ variables $x_1, x_2, \cdots, x_n$  under the constraint $\sum_{i=1}^n x_i=1$ and the theory of Lagrange function and multipliers is often used in evaluating  $\lambda(G)$,  $\lambda (G)$ was called the  Lagrangian of $G$  in several papers \cite{FF, FR84, mubayi06, T}.  Throughout this paper, we also call $\lambda(G)$   the  Lagrangian of $G$.   \end{remark}

The following fact is easily implied by Definition \ref{definitionlagrangian}.

\begin{fact}\label{mono}
Let $G_1$, $G_2$ be $r$-uniform graphs and $G_1\subseteq G_2$. Then $\lambda (G_1) \le \lambda (G_2).$
\end{fact}

In \cite{MS}, Motzkin and Straus provided the following simple expression for the Lagrangian of a 2-graph.

\begin{theorem} {\rm(See \cite{MS}, Theorem 1)} \label{MStheo}
If $G$ is a 2-graph with $n$ vertices in which a largest clique has order $t$ then
$\lambda(G)=\lambda(K^{(2)}_t)={1 \over 2}(1 - {1 \over t})$. Furthermore,  the vector $\vec{x}:=(x_1,x_2,\ldots ,x_n)$ given by $x_{i}:={1 \over t}$ if $i$ is a  vertex in a fixed maximum clique and $x_i=0$ otherwise is an optimal weighting.
\end{theorem}

This result provides a solution to the optimazation problem of this type of homogeneous
quadratic functions over the standard simplex of an Euclidean plane. It is well-known that
Lagrangians of hypergraphs have been proved to be a useful tool in hypergraph extremal
problems, for example, it has been applied in finding Tur\'{a}n densities of hypergraphs in \cite{sidorenko89,FF,mubayi06}. In order to explore the relationship between the Lagrangian of a hypergraph
and the order of its maximum cliques for hypergraphs when the number of edges is in certain ranges, the following two conjectures are proposed in [17].

\begin{conjecture} \label{conjecture1} (See \cite{PZ}, Conjecture 1.3)
Let $m$ and $t$ be positive integers satisfying ${t-1 \choose r} \le m \le {t-1 \choose r} + {t-2 \choose r-1}.$
Let $G$ be an $r$-graph with $m$ edges and  contain a clique of order  $t-1$. Then $\lambda(G)=\lambda([t-1]^{(r)})$.
\end{conjecture}

\begin{conjecture} \label{conjecture2} (See \cite{PZ}, Conjecture 1.4)
Let $m$ and $t$ be positive integers satisfying ${t-1 \choose r} \le m \le {t-1 \choose r} + {t-2 \choose r-1}.$
 Let $G$ be an $r$-graph with $m$ edges and contain no clique of order $t-1$.  Then $\lambda(G) < \lambda([t-1]^{(r)})$.
\end{conjecture}

In \cite{PZ}, we proved that Conjecture \ref{conjecture1} holds for $r=3$.

\begin{theorem} \label{theorem 1} {\rm (See \cite{PZ}, Theorem 1.8)} Let $m$ and $t$ be positive integers satisfying ${t-1 \choose 3} \le m \le {t-1 \choose 3} + {t-2 \choose 2}$. Let $G$ be a $3$-graph with $m$ edges and  contain a clique of order  $t-1$. Then $\lambda(G) = \lambda([t-1]^{(3)}).$
\end{theorem}

For distinct $A, B \in {\mathbb N}^{(r)}$ we say that $A$ is less than $B$ in the {\em colex ordering} iff $max(A \triangle B) \in B$, where \\$A \triangle B:=(A \setminus B)\cup (B \setminus A)$. For example we have $246 < 156$ in ${\mathbb N}^{(3)}$ since $max(\{2,4,6\} \triangle \{1,5,6\}) \in \{1,5,6\}$. In colex ordering, $123<124<134<234<125<135<235<145<245<345<126<136<236<146<246<346<156<256<356<456<127<\cdots .$ Note that the first $t \choose r$ $r$-tuples in the colex ordering of ${\mathbb N}^{(r)}$ are the edges of $[t]^{(r)}$.

Let $C_{r,m}$ denote the $r$-graph with $m$ edges formed by taking the first $m$ sets in the colex ordering of ${\mathbb N}^{(r)}$.
The following  result in \cite{T} states that the value of $\lambda(C_{r,m})$ can be easily figured out when $m$ is in a certain range.

\begin{lemma} {\rm (See \cite{T}, Lemma 2.4 )} \label{LemmaTal7}
For any integers $m,t,$ and $r$ satisfying ${t-1 \choose r} \le m \le {t-1 \choose r} + {t-2 \choose r-1},$
we have $\lambda(C_{r,m}) = \lambda([t-1]^{(r)})$.
\end{lemma}

Note that Conjectures \ref{conjecture1} and \ref{conjecture2} refine the following open conjecture of Frankl and F\"uredi.

\begin{conjecture} (See \cite{FF}, Conjecture 4.1)\label{conjecture} The $r$-graph with $m$ edges formed by taking the first $m$ sets in the colex ordering of ${\mathbb N}^{(r)}$ has the largest Lagrangian of all $r$-graphs with  $m$ edges. In particular, the $r$-graph with $t \choose r$ edges and the largest Lagrangian is $[t]^{(r)}$.
\end{conjecture}

Note that the upper bound ${t-1 \choose r} + {t-2 \choose r-1}$ in Conjecture \ref{conjecture1} is the best possible. For example, if \\$m ={t-1 \choose r}+{t-2 \choose r-1}+1$ then $\lambda(C_{r,m}) > \lambda([t-1]^{(r)})$. To see this,  take  $\vec{x}:=(x_1,\ldots,x_t) \in S$, where \\$x_1=x_2=\cdots=x_{t-2}={1 \over t-1}$ and $x_{t-1}=x_{t}={1 \over 2(t-1)}$,
then $\lambda(C_{r,m})\ge \lambda(C_{r,m},\vec{x} )> \lambda([t-1]^{(r)}).$

In \cite{T}, Talbot proved the following.

\begin{theorem} {\rm(See \cite{T}, Theorem 2.1)} \label{Tal} Let $m$ and $t$ be integers satisfying
${t-1 \choose 3} \le m \le {t-1 \choose 3} + {t-2 \choose 2} - (t-1).$
Then $\lambda(G) \leq \lambda([t-1]^{(3)})$.
\end{theorem}

\begin{theorem} {\rm (See \cite{T}, Theorem 3.1)} \label{TalLemma19}For any $r\geq 4$ there exists constants $\gamma_{r}$ and $k_{0}(r)$ such that if $m$ satisfies $${t-1 \choose r} \le m \le {t-1 \choose r} + {t-2 \choose r-1}-\gamma_{r}(t-1)^{r-2}$$ with $t\geq k_{0}(r)$ and $G$ is an $r$-graph on $t$ vertices with $m$ edges,  then $\lambda(G) \leq \lambda([t-1]^{(r)})$.
\end{theorem}
Note that, Theorems 2.3 and 2.4 in this paper are equivalent to Theorems 2.1 and 3.1 in [18] after shifting $t$ to $t-1$.

Some evidence of Conjectures \ref{conjecture1} and \ref{conjecture2} can be found in \cite{PTZ,PZZ}. In particular,  we proved

\begin{theorem} {\rm(See \cite{PTZ}, Theorem 1.10)}  \label{Corollary 2}
 (a) Let $m$ and $t$ be positive integers satisfying $${t-1 \choose r} \le m \le {t-1 \choose r} + {t-2 \choose r-1}-(2^{r-3}-1)({t-2 \choose r-2}-1).$$  Let $G$ be an $r$-graph on $t$ vertices with $m$ edges and contain a clique of order $t-1.$ Then $\lambda(G) = \lambda([t-1]^{(r)})$.

(b) Let $m$ and $t$ be positive integers satisfying ${t-1 \choose 3} \le m \le {t-1 \choose 3} + {t-2 \choose 2}-(t-2)$.  Let $G$ be a $3$-graph with $m$ edges and without containing a clique of order $t-1$. Then  $\lambda(G) <\lambda([t-1]^{(3)}).$
\end{theorem}
In this paper, we provide upper bounds on the Lagrangian of a 3-graph, a 4-graph, and an r-graph,  respectively,  when the hypergraph contains  dense subgraphs and  the number of edges of the hypergraph is in a certain range. These results support Conjectures \ref{conjecture1}, \ref{conjecture2} and extend Theorem \ref{Tal}.


We will impose one additional condition on any optimal weighting ${\vec x}:=(x_1, x_2, \ldots, x_n)$ for an $r$-graph $G$:
\begin{eqnarray}
 &&|\{i : x_i > 0 \}|{\rm \ is \ minimal, i.e. \ if}  \ \vec y {\rm \ is \ a \ feasible \ weighting \ for \ } G  {\rm \ satisfying }\nonumber \\
 &&|\{i : y_i > 0 \}| < |\{i : x_i > 0 \}|,  {\rm \  then \ } \lambda (G, {\vec y}) < \lambda(G) \label{conditionb}.
\end{eqnarray}

When the theory of Lagrange multipliers is applied to find the optimum of $\lambda(G)$, subject to $\sum_{i=1}^n x_i =1$, note that $\lambda (E_i, {\vec x})$ corresponds to the partial derivative of  $\lambda(G, \vec x)$ with respect to $x_i$. The following lemma gives some necessary conditions of an optimal weighting of  $\lambda(G)$.

\begin{lemma} {\rm (See \cite{FR84}, Theorem 2.1)} \label{LemmaTal5} Let $G:=(V,E)$ be an $r$-graph on the vertex set $[n]$ and ${\vec x}:=(x_1, x_2, \ldots, x_n)$ be an optimal feasible weighting for $G$ with $k$  ($\le n$) non-zero weights $x_1, x_2, \ldots, x_k$ satisfying condition (\ref{conditionb}). Then for every $\{i, j\} \in [k]^{(2)}$,  (a) $\lambda (E_i, {\vec x})=\lambda (E_j, \vec{x})=r\lambda(G)$, (b) there is an edge in $E$ containing both $i$ and $j$.
\end{lemma}
The following definition is also needed.
\begin{definition}
An $r$-graph $G:=(V,E)$  on the vertex set $[n]$ is {\it left-compressed} if $j_1j_2 \cdots j_r \in E$ implies \\$i_1i_2 \cdots i_r \in E$ provided $i_p \le j_p$ for every $p, 1\le p\le r$. Equivalently, an $r$-graph $G:=(V,E)$ is {\it left-compressed} iff $E_{j\setminus i}=\emptyset$ for any $1\le i<j\le n$.
\end{definition}
\begin{remark}\label{r1} (a) In Lemma \ref{LemmaTal5}, part(a) implies that
$$x_j\lambda(E_{ij}, {\vec x})+\lambda (E_{i\setminus j}, {\vec x})=x_i\lambda(E_{ij}, {\vec x})+\lambda (E_{j\setminus i}, {\vec x}).$$
In particular, if $G$ is left-compressed, then
$$(x_i-x_j)\lambda(E_{ij}, {\vec x})=\lambda (E_{i\setminus j}, {\vec x})$$
for any $i, j$ satisfying $1\le i<j\le k$ since $E_{j\setminus i}=\emptyset$.

(b) If  $G$ is left-compressed, then for any $i, j$ satisfying $1\le i<j\le k$,
\begin{equation}\label{enbhd}
x_i-x_j={\lambda (E_{i\setminus j}, {\vec x}) \over \lambda(E_{ij}, {\vec x})}
\end{equation}
holds.  If  $G$ is left-compressed and  $E_{i\setminus j}=\emptyset$ for $i, j$ satisfying $1\le i<j\le k$, then $x_i=x_j$.

(c) By (\ref{enbhd}), if  $G$ is left-compressed, then an optimal feasible weighting  ${\vec x}:=(x_1, x_2, \ldots, x_n)$ for $G$  must satisfy
\begin{equation}\label{conditiona}
x_1 \ge x_2 \ge \ldots \ge x_n \ge 0.
\end{equation}
\end{remark}

In the proofs of our results, we need to consider various left-compressed 3-graphs  on vertex set [t],
 which can be obtained from a Hessian diagram as follows.

A triple  $i_1 i_2 i_3$ is called a {\it descendant  } of a triple  $j_1j_2 j_3$ iff $i_s\le j_s$ for each $1\le s\le 3$, and \\$i_1+i_2+i_3 < j_1+j_2+j_3$. In this case, the triple $j_1j_2j_3$  is called an {\it ancestor } of $i_1 i_2 i_3$.  The triple $i_1 i_2 i_3$   is called a {\it direct descendant} of $j_1j_2 j_3$ if  $i_1 i_2 i_3$ is a descendant of $j_1j_2 j_3$ and $j_1+j_2+j_3=i_1+i_2+i_3 +1$.  We say that $j_1 j_2 j_3$ has lower hierarchy than $i_1i_2 i_3$ if $j_1 j_2 j_3$  is  an ancestor of $i_1i_2i_3$. This is a partial order on the set of all triples.  Fig.1 is a Hessian diagram on all triples on vertex set $[t]$. In this diagram, $i_1 i_2 i_3$ and  $j_1j_2 j_3$ are connected by an edge if $i_1 i_2 i_3$   is  a  direct descendant of $j_1j_2 j_3$.


\begin{figure}
 \centering
  \includegraphics{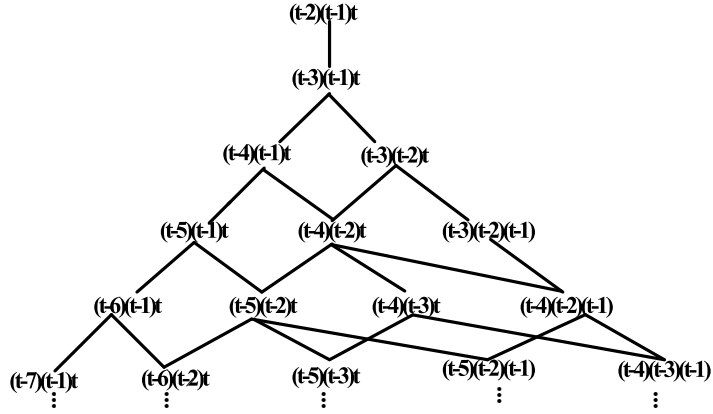}
\caption{Hessian Diagram on $[t]^{(3)}$}
\label{fig:1}       
\end{figure}

\begin{remark}\label{releftcom}
A $3$-graph $G$ is left-compressed iff all descendants of an edge of $G$ are edges of $G$. Equivalently, if a triple is not an edge of $G$, then none of its ancestors will be an edge of $G$.
\end{remark}

\section{The Lagrangians of $3$-graphs Containing Subgraph $K_{t-1}^{(3)-}$}
\label{sec:3}
Let $K_{t-1}^{(3)-}$ denote the hypergraph obtained by $K_{t-1}^{(3)}$ with one edge removed, where  $K_{t-1}^{(3)}$ stands for a complete $3$-graph with $t-1$ vertices.
Denote $\lambda_{(m,t-1)}^{3-}:=\max\{\lambda(G): G $ is a $3$-graph with $m$ edges and $G$
containing  $ K_{t-1}^{(3)-}$ but not containing  $K_{t-1}^{(3)}\}$.
We now prove Theorem \ref{theorem31}.

\begin{theorem} \label{theorem31}
Let $m$ and $t$ be positive integers satisfying ${t-1 \choose 3}\leq m \le {t-1 \choose 3}+{t-2 \choose 2}$. Let $G$ be a $3$-graph with $m$ edges containing $K_{t-1}^{(3)-}$ but not containing $K_{t-1}^{(3)}$.  Then $\lambda(G)<\lambda([t-1]^{(3)})$ for $t\geq 6$.
\end{theorem}

In the proof of Theorem \ref{theorem31},  we need several lemmas.

\begin{lemma}\label{lemma31} Let $m$ and $t$ be positive integers satisfying ${t-1 \choose 3} \le m \le {t-1 \choose 3}+{t-2 \choose 2}.$ Then there exists a left-compressed $3$-graph $G$ with $m$ edges containing $[t-1]^{(3)} \backslash \{(t-3)(t-2)(t-1)\}$ but not containing $[t-1]^{(3)}$ such that $\lambda(G)=\lambda_{(m,t-1)}^{(3)-}$ and there exists an optimal weighting $\vec{x}:=(x_{1},x_{2},\ldots ,x_{n})$ of $G$ satisfying $x_{i}\geq x_{j}$ when $i<j$.
\end{lemma}

The proof of Lemma \ref{lemma31} is similar to the proof of Lemma 3.1 in \cite{PZZ}. However  Lemma 3.1 in \cite{PZZ} cannot be used directly here. For completeness, we give  the proof in Sect. \ref{Lemmas}.

\begin{lemma} {\rm(See \cite{PZZ}, Proposition 3.7 )}\label{lemma32} Let $G$ be a $3$-graph on $t$ vertices with at most ${t-1 \choose 3}+{t-1 \choose 2}$ edges. If $G$ does not contain $K_{t-1}^{(3)}$, then $\lambda(G)<\lambda([t-1]^{(3)})$ for $6\leq t\leq 12$.
\end{lemma}

\begin{lemma}\label{lemma33}  Let $G$ be a left-compressed 3-graph containing $[t-1]^{(3)}\backslash\{(t-3)(t-2)(t-1)\}$ but not containing $[t-1]^{(3)}$ with $m$ edges such that  $\lambda(G)=\lambda_{(m,t-1)}^{3-}$. Let $\vec{x}:=(x_{1},x_{2},\ldots ,x_{n})$ be an optimal weighting of $G$ and $k$ be the number of positive weights in $\vec{x}$, then $\lambda(G)<\lambda([t-1])^{(3)}$ or  $|[k-1]^{(3)}\backslash E|\leq k-2.$
\end{lemma}
The proof of Lemma \ref{lemma33} is similar to  Lemma 3.2 in \cite{PZZ}. However  Lemma 3.2 in \cite{PZZ} cannot be used directly here. For completeness, we give the details of the proof in Sect. \ref{Lemmas}.

\bigskip
\noindent {\em Proof of Theorem \ref{theorem31}} Let $m$ and $t$ be positive integers satisfying ${t-1 \choose 3} \le m \le {t-1 \choose 3}+{t-2 \choose 2}.$  Let $G:=(V,E)$ be a $3$-graph with $m$ edges containing $K_{t-1}^{(3)-}$ but not containing $K_{t-1}^{(3)}$ such that $\lambda(G)=\lambda_{(m,t-1)}^{3-}$. Let \\$\vec{x}:=(x_{1},x_{2},\ldots ,x_{n})$ be an optimal weighting of $G$ and $k$ be the number of non-zero weights in $\vec{x}$. By Lemma \ref{lemma31}, we can assume that $G$ is left-compressed and contains $[t-1]^{(3)}\backslash\{(t-3)(t-2)(t-1)\}$ but not contain $[t-1]^{(3)}$ and $x_1 \ge x_2 \ge \ldots \ge x_k >x_{k+1}=\ldots=x_{n}=0$. Since $\vec{x}$ has only $k$ positive weights, we can assume that $G$ is on $[k]$.

Now we proceed to show that $\lambda(G) <\lambda([t-1]^{(3)})$. By Lemma \ref{lemma32}, Theorem \ref{theorem31} holds when $t\leq 12$. Next we assume $t\geq 13$. If  $\lambda(G)\ge \lambda([t-1]^{(3)})$, then $k\geq t$. Otherwise $k\leq t-1$, since $G$ does not contain $[t-1]^{(3)}$, then  $\lambda(G)<\lambda([t-1]^{(3)})$.

Since $G$ is left-compressed and $1(k-1)k\in E$, then $\vert[k-2]^{(2)}\cap E_{k}\vert\geq 1$.
If $k\ge t+1$, then applying Lemma \ref{lemma33},  we have
\begin{eqnarray}\label{eqs31}
m=\vert E\vert &=&\vert E\cap [k-1]^{(3)}\vert +\vert [k-2]^{(2)}\cap E_k \vert +\vert E_{(k-1)k}\vert \nonumber\\
&\ge & {t \choose 3}-(t-1) +2 \nonumber\\
&\ge& {t-1 \choose 3} + {t-2 \choose 2}+1,
\end{eqnarray}
which contradicts to the assumption that $m\le {t-1 \choose 3} + {t-2 \choose 2}$. Recall that $k\ge t$, so we have
$$k=t.$$
Since $\lambda_{(m,t-1)}^{3-}$ does not decrease as $m$ increases, it is sufficient to show  the case that $m={t-1 \choose 3} + {t-1 \choose 2}$.
Let $G':=G\bigcup \{(t-3)(t-2)(t-1)\} \backslash \{1(t-1)t\}$. If we can prove that $\lambda(G,\vec{x})<\lambda(G',\vec{x})$, then since $G'$ contains $[t-1]^{(3)}$ and $G'$ has ${t-1 \choose 3} + {t-1 \choose 2}$ edges, we have $\lambda(G',\vec{x})\leq \lambda (G')=\lambda([t-1]^{(3)})$. Consequently, \\$\lambda(G)<\lambda([t-1]^{3})$. Now we show that $\lambda(G,\vec{x})<\lambda(G',\vec{x})$. Note that
\begin{equation}\label{compare}
\lambda(G',\vec{x})- \lambda(G,\vec{x})=x_{t-3}x_{t-2}x_{t-1}-x_{1}x_{t-1}x_{t}.
\end{equation}

By Remark \ref{r1}(b), we have
 \begin{equation}\label{eqs322}
 x_{1}=x_{t-3}+\frac{\lambda(E_{1\setminus(t-3)},\vec{x})}{\lambda(E_{1(t-3)},\vec{x})},
\end{equation}
and
\begin{equation}\label{eqs333}
 x_{t-2}=x_{t}+\frac{\lambda(E_{(t-2)\backslash t},\vec{x})}{\lambda(E_{(t-2) t},\vec{x})}.
\end{equation}

Combining equations (\ref{compare}), (\ref{eqs322}) and (\ref{eqs333}), we get
\begin{eqnarray}\label{compare2}
\lambda(G',\vec{x})- \lambda(G,\vec{x})&=&x_{t-3}(x_{t}+\frac{\lambda(E_{(t-2)\backslash t},\vec{x})}{\lambda(E_{(t-2) t},\vec{x})})x_{t-1}-(x_{t-3}+\frac{\lambda(E_{1\setminus(t-3)},\vec{x})}{\lambda(E_{1(t-3)},\vec{x})})x_{t-1}x_{t} \nonumber \\
&=&\frac{\lambda(E_{(t-2)\backslash t},\vec{x})}{\lambda(E_{(t-2) t},\vec{x})}x_{t-3}x_{t-1}-\frac{\lambda(E_{1\setminus(t-3)},\vec{x})}{\lambda(E_{1(t-3)},\vec{x})}x_{t-1}x_{t}.
\end{eqnarray}
By Remark \ref{r1}(b)
 \begin{eqnarray}
x_{1}=x_{t-2}+\frac{\lambda(E_{1\backslash (t-2)},\vec{x})}{\lambda(E_{1(t-2)},\vec{x})}\leq x_{t-2}+\frac{x_{t-3}x_{t-1}+(x_{2}+\cdots +x_{t-3})x_{t}}{x_{2}+ \cdots+x_{t-3}+x_{t}}<x_{t-2}+x_{t-1}+x_{t}.
\end{eqnarray}
Hence $\lambda(E_{1(t-3)},\vec{x})-\lambda(E_{(t-2)t},\vec{x})\ge x_{t-2}+x_{t-1}+x_{t}-x_{1}>0$ and $\lambda(E_{1(t-3)},\vec{x})>\lambda(E_{(t-2)t},\vec{x})$. Clearly \\$x_{t-3}>x_{t}$ since $(t-5)(t-1)\in E_{(t-3)\backslash t}$. Therefore to show that $\lambda(G,\vec{x})<\lambda(G',\vec{x})$, it is sufficient to show that
\begin{equation}\label{compare3}
\lambda(E_{(t-2)\backslash t},\vec{x})\ge \lambda(E_{1\setminus(t-3)},\vec{x}).
\end{equation}

If $(t-6)(t-1)t\in E$, then all triples in $[t]^{(3)}\setminus\{(t-3)(t-2)(t-1), ijt, {\rm where} \ t-5\le i<j\le t-1\}$ are edges in $G$ since $G$ is left-compressed. If $E\neq [t]^{(3)}\setminus\{(t-3)(t-2)(t-1), ijt, {\rm where} \ t-5\le i<j\le t-1\}$, then
$m>{t \choose 3}-11\ge {t-1 \choose 3}+{t-2 \choose 2}$(recall that $t\geq 13$.), which is a contradiction. Therefore, either\\ $E=[t]^{(3)}\setminus\{(t-3)(t-2)(t-1), ijt, {\rm where} \ t-5\le i<j\le t-1\}$ or $(t-6)(t-1)t\notin E$.

If $E=[t]^{(3)}\setminus\{(t-3)(t-2)(t-1), ijt, {\rm where} \ t-5\le i<j\le t-1\}$, then
$$\lambda(E_{(t-2)\backslash t},\vec{x})=x_{t-5}x_{t-1}+x_{t-5}x_{t-3}+x_{t-5}x_{t-4}+x_{t-4}x_{t-3}+x_{t-4}x_{t-1},$$
and
$$\lambda(E_{1\setminus(t-3)},\vec{x})=x_{t-2}x_{t-1}+x_{t-5}x_{t}+x_{t-4}x_{t}+x_{t-2}x_{t}+x_{t-1}x_{t}.$$
Clearly (\ref{compare3}) holds in this case.

If $(t-6)(t-1)t\notin E$, then
\begin{eqnarray*}
\lambda(E_{(t-2)\backslash t},\vec{x})&\ge&x_{t-3}\lambda(E_{(t-3)(t-2)}\cap E_{(t-3)t}^{c},\vec{x})+x_{t-4}x_{t-1}+x_{t-5}x_{t-1}+x_{t-6}x_{t-1}\nonumber\\
 &= & x_{t-3}\lambda(E_{(t-3)t}^{c},\vec{x})+x_{t-4}x_{t-1}+x_{t-5}x_{t-1}+x_{t-6}x_{t-1}-x_{t-3}x_{t-2}-x_{t-3}x_{t-1}\nonumber\\
 &\geq & x_{t-3}\lambda(E_{(t-3)t}^{c},\vec{x})+x_{t-5}x_{t-1}+x_{t-6}x_{t-1}-x_{t-3}x_{t-2}\nonumber\\
 &=&x_{t-3}(\lambda(E_{(t-3)t}^{c},\vec{x})-x_{t-2})+x_{t-5}x_{t-1}+x_{t-6}x_{t-1},
 \end{eqnarray*}
and
\begin{eqnarray*}
 \lambda(E_{1\setminus(t-3)},\vec{x})&=&x_{t}\lambda(E_{(t-3)t}^{c},\vec{x})+x_{t-2}x_{t-1}\nonumber\\
 &=&x_{t}(\lambda(E_{(t-3)t}^{c},\vec{x})-x_{t-2})+x_{t-2}x_{t-1}+x_{t-2}x_{t}.
 \end{eqnarray*}
Clearly (\ref{compare3}) holds in this case.

This completes the proof of Theorem \ref{theorem31}.\qed

\begin{remark} Note that for $t\leq 5$, the left-compressed $3$-graph with ${t-1 \choose 3}+{t-1 \choose 2}$ edges always contains $K_{t-1}^{(3)}$. Combining Theorems \ref{Corollary 2} and \ref{theorem31}, we have that, if $G$ is a $3$-graph containing $K_{t-1}^{(3)-}$ with at most ${t-1 \choose 3}+{t-1 \choose 2}$ edges, then $\lambda(G)\leq \lambda([t-1]^{(3)})$.
\end{remark}

Also, applying  Theorem \ref{theorem31},  we derive two easy corollaries that support Conjecture 2.2.

\begin{corollary}\label{tvertices} Let $m$ and $t$ be positive integers satisfying ${t-1 \choose 3} \le m \le {t-1 \choose 3} + {t-2 \choose 2}$.  Let $G:=(V,E)$ be a left-compressed 3-graph on the vertex set [t] with $m$ edges and not containing a clique of order $t-1$. If  $|E_{(t-1)t}|\le3$, then $\lambda (G)<\lambda ([t-1]^{(3)})$.
\end{corollary}

\noindent{\em Proof } Because $\lambda_{(m,t-1)}^{3-}$ doesn't decrease as $m$ increases, we can assume that $m={t-1 \choose 3} + {t-2 \choose 2}$.
Since \\$G:=(V, E)$ does not contain $[t-1]^{(3)}$ and $G$ is left-compressed, then $(t-3)(t-2)(t-1)\notin E$. If $|E_{(t-1)t}|=1$, then $G$ must contain $[t-1]^{(3)}$. Therefore, $|E_{(t-1)t}|=2$ or $3$.

If $t\leq 5$, Theorem \ref{tvertices} clearly holds.  Next, we assume $t\geq 6$ and distinguish two cases.\\
Case 1. $|E_{(t-1)t}|=2$.
Note that $G$ is left-compressed, in view of Fig.1, $$E=[t]^{(3)}\setminus\{3(t-1)t,4(t-1)t,\cdots (t-2)(t-1)t,(t-3)(t-2)(t-1),(t-3)(t-2)t\}.$$
Case 2. $|E_{(t-1)t}|=3$.
In this case, since $G$ is left-compressed, in view of Fig.1, we only need to consider $E=[t]^{(3)}\setminus\{4(t-1)t,\cdots (t-2)(t-1)t,(t-3)(t-2)(t-1),(t-3)(t-2)t,(t-4)(t-2)t\}$.

In both cases, left-compressed 3-graph G does not contain the edge (t-3)(t-2)(t-1). Thus, the conditions in Theorem 3.1 are satisfied. Therefore, we are done.
 \qed

The next corollary states that if 3-graph $G$ contains a dense subgraph close to the structure in $C_{3,m}$, then we have $\lambda(G) < \lambda([t-1]^{(3)})$.

\begin{corollary}\label{symmetirc_difference} Let $m$ and $t$ be positive integers satisfying ${t-1 \choose 3} \le m \le {t-1 \choose 3} + {t-2 \choose 2}$.  Let $G:=(V,E)$ be a left-compressed 3-graph on the vertex set [t] with $m$ edges and not containing a clique of size $t-1$, and $|E(G) \Delta E(C_{3,m})| \le 6$. Then, $\lambda (G) < \lambda ([t-1]^{(3)})$.
\end{corollary}

\noindent{\em Proof }  If  $m\le {t-1 \choose 3}+{t-2 \choose 2}$, then  $|E_{(t-1)t}|\leq 3$, since otherwise $|E(G) \Delta E(C_{3,m})| >6$. Applying  Corollary \ref{symmetirc_difference}, we have $\lambda(G) < \lambda([t-1]^{(3)})$. \qed


\section{The Lagrangians of Hypergraphs Containing A Clique of Order $t-2$ or $t-1$}
\label{sec:4}
In this section, we prove the following.

\begin{theorem} \label{theorem41} Let $m$ and $t$ be positive integers satisfying ${t-1 \choose 3} \le m \le {t-1 \choose 3} + {t-2 \choose 2}-\frac{t-2}{2}$. Let $G$ be a $3$-graph with $m$ edges and $G$ contain the maximum clique of order  $t-2$. Then $\lambda(G)<\lambda([t-1]^{(3)})$.
\end{theorem}

\begin{theorem} \label{theorem42} Let $m$ and $t$ be positive integers satisfying ${t-1 \choose r} \le m \le {t-1 \choose r} + {t-2 \choose r-1}-2^{r-2}({t-2 \choose r-2}-1)$.  Let $G$ be an $r$-graph on $t$ vertices with $m$ edges and with a clique of order $t-2$. Then $\lambda(G) \leq \lambda([t-1]^{(r)})$.
\end{theorem}

\begin{theorem} \label{theorem43} Let $m$ and $t$ be positive integers satisfying ${t-1 \choose 4} \le m \le {t-1 \choose 4} + {\lfloor \frac{t-2}{2} \rfloor \choose 3}$.  Let $G$ be a $4$-graph with $m$ edges and a clique of order $t-1$. Then $\lambda(G) = \lambda([t-1]^{(4)})$.
\end{theorem}
Here ${t-1 \choose 4} + {\lfloor \frac{t-2}{2} \rfloor \choose 3}$ is not the best upper bound that we can obtain. This bound is for simplicity of the proof.

Denote $\lambda_{(m,p)}^{r}:=\max\{\lambda(G): G {\rm \ is \ an \ } r{\rm -graph  \ with \ } m$
edges and $G$ contains  a  maximum  clique  of
order $p \}.$

Similar to the proof of Lemma 3.1 in \cite{PZZ}, we can prove the following lemma. We will give the proof in Sect. \ref{Lemmas}.
\begin{lemma}\label{lemma41} Let $m$ and $t$ be positive integers satisfying $${t-1 \choose 3} \le m \le {t-1 \choose 3}+{t-2 \choose 2}-{t-2 \over 2}.$$
Then there exists a left-compressed $3$-graph $G$ with $m$ edges containing the maximum clique $[t-2]^{(3)}$ such that $\lambda(G)=\lambda_{(m,t-2)}^{3}$.
\end{lemma}

Similar to the proof of Lemma 3.2 in \cite{PZZ}, we have the following lemma. For completeness, we will give the proof in Sect. \ref{Lemmas}.

\begin{lemma}\label{lemma42}  Let $G$ be a left-compressed 3-graph containing the maximum clique $[t-2]^{(3)}$ with $m$ edges such that  $\lambda(G)=\lambda_{(m,t-2)}^{3}$. Let $\vec{x}:=(x_{1},x_{2},\ldots ,x_{n})$ be an optimal weighting of $G$ and $k$ be the number of positive weights in $\vec{x}$, then $\lambda(G)<\lambda([t-1])^{(3)}$ or  $|[k-1]^{(3)}\backslash E|\leq k-2.$
\end{lemma}

We also need the following lemma whose proof is similar to  Lemma 2.7 in \cite{T} and Lemma 3.3 in \cite{PZ}.
We will  give it in Sect. \ref{Lemmas}.

\begin{lemma}\label{lemma43} Let $m$ and $t$ be positive integers satisfying ${t-1 \choose 3} \le m \le {t-1 \choose 3} + {t-2 \choose 2}-\frac{t-2}{2}$. Let $G$ be a left-compressed 3-graph on the vertex set $[t]$ and contain the maximum clique $[t-2]^{(3)}$ with $m$ edges such that  $\lambda(G)=\lambda_{(m,t-2)}^{3}$. Assume $b:=|E_{(t-1)t}|$, then $\lambda(G)< \lambda ([t-1]^{(3)})$ or $$|[t-2]^{(2)}\backslash E_{t}|\leq b.$$
\end{lemma}

\noindent {\em Proof of Theorem \ref{theorem41}}
Let $m$ and $t$ be positive integers satisfying ${t-1 \choose 3} \le m \le {t-1 \choose 3}+{t-2 \choose 2}-\frac{t-2}{2}.$  Clearly we can assume that $t\geq 5$. Let $G:=(V,E)$ be a $3$-graph with $m$ edges containing a maximum clique  of order $t-2$ such that $\lambda(G)=\lambda_{(m,t-2)}^{3}$. Let $\vec{x}:=(x_{1},x_{2},\ldots ,x_{n})$ be an optimal weighting of $G$ and $k$ be the number of non-zero weights in $\vec{x}$. By Lemma \ref{lemma41}, we can assume that $G$ is left-compressed with the maximum clique $[t-2]^{(3)}$ and $x_1 \ge x_2 \ge \ldots \ge x_k >x_{k+1}=\ldots=x_{n}=0$. Since $\vec{x}$ has only $k$ positive weights, we can assume that $G$ is on $[k]$.

Now we proceed to show that $\lambda(G) <\lambda([t-1]^{(3)})$.  If  $\lambda(G)\ge \lambda([t-1]^{(3)})$, then $k\geq t$. Otherwise $k\leq t-1$, since $G$ does not contain $[t-1]^{(3)}$, then  $\lambda(G)<\lambda([t-1]^{(3)})$. By Lemma \ref{LemmaTal5}(a), $k-1$ and $k$ appear in some common edge $e \in E$. Recall that $E$ is left-compressed, so $1(k-1)k\in E$. Define $b:=max\{i:i(k-1)k\in E\}$. Because $E$ is left-compressed, $E_{i\setminus j}=\emptyset$ for $1\le i<j\le b$. Hence, by Remark \ref{r1}(a), we have $x_1=x_2=\cdots=x_b$. Clearly, $b\leq k-5$.

Since $G$ is left-compressed and $1(k-1)k\in E$, then $\vert[k-2]^{(2)}\cap E_{k}\vert\geq 1$.
So applying Lemma \ref{lemma42}, similar to  (\ref{eqs31}), we have $k=t$.

Since $k=t$, we can assume that $G$ is on $[t]$. 
By Remark \ref{r1}(b), we have $$x_{1}= x_{t-3}+\frac{\lambda(E_{1\setminus(t-3)},\vec{x})}{\lambda(E_{1(t-3)},\vec{x})}.$$
Recall that $G$ contains a clique order of $t-2$, we have $$\lambda(E_{1\setminus(t-3)},\vec{x})=x_{t-1}\lambda(E_{(t-3)(t-1)}^{c},\vec{x})+x_{t}\lambda(E_{(t-3)t}^{c},\vec{x})-x_{t-1}x_{t}.$$
Hence
$$
x_{1}< x_{t-3}+\frac{x_{t-1}\lambda(E_{(t-3)(t-1)}^{c},\vec{x})+x_{t}\lambda(E_{(t-3)t}^{c},\vec{x})}{\lambda(E_{1(t-3)},\vec{x})}.
$$

Since for $i\neq t-1, t-2, t-3$, $i \in E_{(t-3)t}^{c}$ implies that $i(t-3) \in [t-2]^{(2)}\backslash E_{t}$ and $i \in E_{(t-2)t}^{c}$ implies that $i(t-2) \in [t-2]^{(2)}\backslash E_{t}$, $t-1 \in E_{(t-3)t}^{c}$, $t-1 \in E_{(t-2)t}^{c}$,
and \\$t-2\in E_{(t-3)t}^{c}$, $t-3\in E_{(t-2)t}^{c}$ and $(t-2)(t-3)\in [t-2]^{(2)}\backslash E_{t}$, applying Lemma \ref{lemma43}, then
$$|E_{(t-3)t}^{c}|+ |E_{(t-2)t}^{c}|\leq |[t-2]^{(2)}\backslash E_{t}|+3\leq b+3.$$
Note that $b\leq t-5$ and $$|E_{(t-3)t}^{c}|\leq |E_{(t-2)t}^{c}|,$$
So
$$|E_{(t-3)t}^{c}|\leq\frac{b+3}{2}\leq \frac{t-2}{2}.$$
Since $G$ is left-compressed, then
$$|E_{(t-3)(t-1)}^{c}|\leq|E_{(t-3)t}^{c}|\leq \frac{t-2}{2}.$$
So
\begin{eqnarray*}
x_{1}&<& x_{t-3}+\frac{x_{t-1}\lambda(E_{(t-3)(t-1)}^{c},\vec{x})+x_{t}\lambda(E_{(t-3)t}^{c},\vec{x})}{\lambda(E_{1(t-3)},\vec{x})}\\
&\leq& x_{t-3}+\frac{x_{t-1}\frac{t-2}{2}\frac{\lambda(E_{1(t-3)},\vec{x})}{t-2}+x_{t}\frac{t-2}{2}\frac{\lambda(E_{1(t-3)},\vec{x})}{t-2}}{\lambda(E_{1(t-3)},\vec{x})}\\
&\leq& 2x_{t-3}.
\end{eqnarray*}
This implies
$$2x_{t-3}x_{t-2}x_{t-1}-x_{1}x_{t-1}x_{t}>0.$$
Let$C:=[t-1]^{(3)}\backslash E$ be all triples  containing $t-1$ not in $E$,

$E':=E\bigcup C\backslash \{(b-\lfloor\frac{|C|}{2}\rfloor+1)(t-1)t,(b-\lfloor\frac{|C|}{2}\rfloor+2)(t-1)t,\ldots,b(t-1)t\}$and $G':=([t]^{(3)},E')$. Then
\begin{eqnarray*}
\lambda(G',x)-\lambda(G,x)&=& \lambda(C,x)-\lfloor\frac{|C|}{2}\rfloor x_{1}x_{t-1}x_{t}\\
&\geq& |C|x_{t-3}x_{t-2}x_{t-1}-\lfloor\frac{|C|}{2}\rfloor x_{1}x_{t-1}x_{t}\\
&\geq & \frac{|C|}{2}(2x_{t-3}x_{t-2}x_{t-1}-x_{1}x_{t-1}x_{t})>0.
\end{eqnarray*}
So $\lambda(G,x)<\lambda(G',x)$. Because
\begin{eqnarray*}
|E'|&=&|E|+|C|-\lfloor\frac{|C|}{2}\rfloor\leq |E|+\frac{|C|}{2}+1\\
&\leq& {t-1 \choose 3} + {t-2 \choose 2}-\frac{t-2}{2}+\frac{t-4}{2}+1\\
&=& {t-1 \choose 3} + {t-2 \choose 2}.
\end{eqnarray*}
and $G'$ contains a clique of order  $t-1$, we have $\lambda(G',x)\leq \lambda(G')=\lambda([t-1]^{(3)})$ by Theorem \ref{theorem 1}. Hence
$\lambda(G,x)<\lambda(G',x)\leq \lambda([t-1]^{(3)})$. This proves Theorem \ref{theorem41}. \qed




The following lemma implies that we only need to consider left-compressed $r$-graphs when Theorem \ref{theorem42} is proved.
The proof is given in Sect. \ref{Lemmas}.

\begin{lemma}\label{lefcom} Let $m$ and $t$ be positive integers satisfying $${t-1 \choose r} \le m \le {t \choose r}-1.$$Then there exists a left-compressed $G$ with $m$ edges containing the  clique $[t-2]^{(r)}$ such that $\lambda(G)=\lambda_{(m,t-2)}^{r}$ and there exists an optimal weighting $\vec{x}:=(x_{1},x_{2},\ldots ,x_{n})$ of $G$ satisfying $x_{i}\geq x_{j}$ when $i<j$.
\end{lemma}
We also need the following in the proof of Theorem \ref{theorem42} and Theorem \ref{theorem43}
\begin{lemma} {\rm (See \cite{PTZ}, Theorem 3.4)}\label{theoremptz} Let $r\ge 3$ and $t\ge r+2$ be positive integers.
Let $G$ be a left-compressed $r$-graph  on $t$ vertices satisfying $|[t-2]^{(r-1)} \backslash E_t| \ge 2^{r-3}|E_{(t-1)t}|$. Then

(a) If $G$ contains $[t-1]^{(r)}$, then  $\lambda(G) = \lambda([t-1]^{(r)})$,

(b) If $G$  does not contain $[t-1]^{(r)}$, then  $\lambda(G) <\lambda([t-1]^{(r)})$.
\end{lemma}
\noindent {\em Proof of Theorem \ref{theorem42}}
 Let $m$ and $t$ be positive integers satisfying ${t-1 \choose r} \le m \le {t-1 \choose r} + {t-2 \choose r-1}-2^{r-2}({t-2 \choose r-2}-1)$.  Let $G$ be an $r$-graph with $m$ edges and $t$ vertices with a clique order of $t-2$. By Lemma \ref{lefcom} we can assume $G$ is left-compressed.
 By Lemma \ref{theoremptz}, it is sufficient to show that $|[t-2]^{(r-1)} \backslash E_t| \ge 2^{r-3}|E_{(t-1)t}|$.  If not, then  $|[t-2]^{(r-1)} \backslash E_t| < 2^{r-3}|E_{(t-1)t}|$ and $|[t-2]^{(r-1)} \backslash E_{t-1}| \leq|[t-2]^{(r-1)} \backslash E_t| < 2^{r-3}|E_{(t-1)t}|$. Since $G$  contains the clique $[t-2]^{(r)}$, then

\begin{eqnarray*}
m&=&{t-2 \choose r}+2{t-2 \choose r-1}-|[t-2]^{(r-1)} \backslash E_t|-|[t-2]^{(r-1)} \backslash E_{t-1}|+|E_{(t-1)t}|\\
&>&{t-1 \choose r}+{t-2 \choose r-1}-{t-2 \choose r-2}-(2^{r-2}-1)|E_{(t-1)t}|+1 \\
&\ge &{t-1 \choose r}+{t-2 \choose r-1}-2^{r-2}({t-2 \choose r-2}-1).
\end{eqnarray*}
since $|E_{(t-1)t}|\le {t-2 \choose r-2}-1$, this is a contradiction. Note that, if $|E_{(t-1)t}|= {t-2 \choose r-2}$, then $E=[t]^{(r)}$ since $G$ is left-compressed and $m={t \choose r}$, which results in a contradiction too. This proves Theorem \ref{theorem42}.\qed

\begin{remark} \label{remark}
Lemma \ref{theoremptz}(b) and Theorem \ref{theorem42} imply that if $m$ and $t$ are positive integers satisfying $${t-1 \choose r} \le m \le {t-1 \choose r} + {t-2 \choose r-1}-2^{r-2}({t-2 \choose r-2}-1)$$ and $G$ is a $r$-graph on $t$ vertices with $m$ edges and with a maximum clique of order $t-2$. Then $$\lambda(G)<\lambda([t-1]^{(r)}).$$
\end{remark}
\noindent {\em Proof of Theorem \ref{theorem43}}
 Let $m$ and $t$ be positive integers satisfying ${t-1 \choose 4} \le m \le {t-1 \choose 4} + {\lfloor \frac{t-2}{2} \rfloor \choose 3}$.  Let $G$ be a $4$-graph with $m$ edges and a clique of order $t-1$.  Since it contains a clique of order t-1, without loss of generality, we may assume that it contains $[t-1]^{(4)}$. Since $G$ contains $[t-1]^{(4)}$, we have $\lambda(G)\geq \lambda([t-1]^{(4)}).$ Next we prove that $\lambda(G)\leq \lambda([t-1]^{(4)}).$

  Let $\vec{x}:=(x_{1},x_{2},\ldots ,x_{n})$ be an optimal weighting of $G$ and $k$ be the number of non-zero weights in $\vec{x}$. If $k\leq t-1$, clearly $\lambda(G)\leq \lambda([t-1]^{(4)}).$ Assume that $k\geq t$.  Recall that ${t-1 \choose 4} \le m \le {t-1 \choose 4} + {\lfloor \frac{t-2}{2} \rfloor \choose 3}$ and $G$ contains $[t-1]^{(4)}$, hence $|E_k|\leq  {\lfloor \frac{t-2}{2} \rfloor \choose 3}$.
 By Fact \ref{mono}, Lemma \ref{LemmaTal5}(a) and Theorem \ref{Tal}, we have
 $$\lambda(G,\vec{x})=\frac{1}{4}\lambda(E_k,\vec{x})\leq \frac{1}{4}{\lfloor \frac{t-2}{2} \rfloor \choose 3}(\frac{1}{\lfloor\frac{t-2}{2}\rfloor})^{3}\leq \frac{(t-4)(t-6)}{24(t-2)^2}<\frac{(t-2)(t-3)(t-4)}{24(t-1)^3}=\lambda([t-1]^{(4)}).$$\qed

\begin{remark}
Also, note that Theorem \ref{theorem31}, Theorem \ref{theorem41}, and Remark \ref{remark} provide further evidence for Conjecture \ref{conjecture2}. Theorem \ref{theorem43} provide further evidence for Conjecture \ref{conjecture1}.
\end{remark}
\section{Remarks}
\label{sec:5}
Frankl and F\"uredi \cite{FF} asked the following question: Given $r \ge 3$ and $m \in {\mathbb N}$ how large can the Lagrangian of an $r$-graph with $m$ edges be?
Conjecture \ref{conjecture} proposes a  solution to the question mentioned above.

Denote
\begin{eqnarray*}
&&\lambda_m^r:=\max\{\lambda(G): G {\rm \ is \ an \ } r-{\rm graph \ with \ } m {\rm \ edges }\}.\\
\end{eqnarray*}
The following lemma implies that we only need to consider left-compressed $r$-graphs when Conjecture \ref{conjecture} is explored.
\begin{lemma} {\rm (See \cite{T}, Lemma 2.3)}\label{lemmaleftcompress} There exists a left-compressed $r$-graph
$G$ with $m$ edges  such that $$\lambda(G)=\lambda_m^r.$$
\end{lemma}
We extend Theorem \ref{Tal} in Theorem \ref{theorem51} which is a corollary of Theorem \ref{theorem31}.
\begin{theorem} \label{theorem51} Let $m$ and $t$ be positive integers satisfying ${t-1 \choose 3} \le m \le {t-1 \choose 3} + {t-2 \choose 2}-(t-4)$. Then Conjecture \ref{conjecture} is true for $r=3$ and this value of $m$.
\end{theorem}

\noindent {\em Proof } Let $\vec{x}:=(x_{1},x_{2},\ldots ,x_{n})$ be an optimal weighting for $G$ and $k$ be the number of positive weights in $\vec{x}$. We can assume that $G$ is left-compressed by Lemma \ref{lemmaleftcompress}. So $x_1 \ge x_2 \ge \ldots \ge x_k >x_{k+1}=\ldots=x_{n}=0$ by Remark \ref{r1}(c). Since $\vec{x}$ has only $k$ positive weights, we can assume that $G$ is on vertex set $[k]$.

Now we proceed to show that $\lambda(G) \le \lambda([t-1]^{(3)})$.  If  $\lambda(G)>\lambda([t-1]^{(3)})$, then $k\geq t$  since otherwise $k\leq t-1$ and then $\lambda(G)\leq\lambda([t-1]^{(3)})$. Next we apply the following lemma.
 \begin{lemma} {\rm (See \cite{T}, Lemma 2.5)} \label{LemmaTa16}.
Let $m$ be a positive integer. Let $G$ be a left-compressed 3-graph with $m$ edges such that  $\lambda(G)=\lambda_{m}^{3}$. Let $\vec{x}:=(x_{1},x_{2},\ldots ,x_{n})$ be an optimal weighting for $G$ and $k$ be the number of non-zero weights in $\vec{x}$, then $$|[k-1]^{(3)}\backslash E|\leq k-2.$$
\end{lemma}
So similar to  (\ref{eqs31}), we have $k=t$. Next we need the following lemma whose proof follows the lines of Lemma 2.5 in \cite{T}. For completeness, we give  the proof in Sect. \ref{Lemmas}.
\begin{lemma}\label{Lemma51} Let $G$ be a left-compressed 3-graph on the vertex set $[t]$ with $m$ edges where $${t-1 \choose 3} \le m \le {t-1 \choose 3} + {t-2 \choose 2},$$ and $\lambda(G)=\lambda_{m}^{3}$. Let  $\vec{x}:=(x_{1},x_{2},\ldots ,x_{t})$ be an optimal weighting for $G$. Then

$$|[t-1]^{(3)}\backslash E|\leq t-3, \rm or \ \lambda(G)\leq \lambda ([t-1]^{(3)}).$$
\end{lemma}
Assume Lemma \ref{Lemma51} holds, we continue the proof of Theorem \ref{theorem51}. If $\lambda(G)>\lambda ([t-1]^{(3)})$, then \\ $|[t-1]^{(3)}\backslash E|\leq t-3$ by Lemma \ref{Lemma51}, we add any $|[t-1]^{(3)}\setminus E|-1$ triples in $[t-1]^{(3)}\setminus E$ to $E$ and let the new $3$-graph be $G'$. Then $G'$ contains $K_{t-1}^{(3)-}$, the number of edges in $G'$ is at most
${t-1 \choose 3} + {t-2 \choose 2}$ and $\lambda(G')\geq \lambda(G)$. Applying Theorem \ref{Corollary 2} and Theorem \ref{theorem31}, $\lambda(G')\leq\lambda([t-1]^{(3)})$. Therefore $\lambda(G)\leq\lambda([t-1]^{(3)})=\lambda(C_{3,m})$ by Lemma \ref{LemmaTal7}. This completes the proof of Theorem \ref{theorem51}.\qed

\section{\bf Proofs of Some Lemmas} \label{Lemmas}
Proof techniques of lemmas in this section follow from proof techniques of some lemmas in \cite{T,PTZ,PZZ}. As mentioned earlier, lemmas in those papers cannot be applied directly to situations in this paper. For completeness, we give the proof of these lemmas in this section.

\noindent{\em Proof of Lemma \ref{lemma31}} Let $G$ be a 3-graph on the vertex set $[n]$ with $m$ edges containing $K_{t-1}^{(3)-}$  but not containing $K_{t-1}^{(3)}$ such that $\lambda(G)=\lambda_{(m,t-1)}^{3-}$. We call such a 3-graph $G$ an extremal 3-graph for $m$ and $t-1$. Let \\${\vec x}:=(x_1, x_2, \ldots, x_n)$ be an optimal weighting of $G$. We can assume that $x_i\ge x_j$ when $i<j$ since otherwise we can just relabel the vertices of $G$ and obtain another  extremal $3$-graph for $m$ and $t-1$ with an optimal weighting ${\vec x}:=(x_1, x_2, \ldots, x_n)$ satisfying $x_i\ge x_j$ when $i<j$. Next we obtain a new $3$-graph $G'$ from $G$ by performing the following:

\begin{enumerate}

\item If $(t-3)(t-2)(t-1) \in E(G)$, then there is at least one triple  in $[t-1]^{(3)}\setminus E(G)$,  we replace\\ $(t-3)(t-2)(t-1)$ by this triple;

\item If  an  edge in $G$ has  a  descendant other than $(t-3)(t-2)(t-1)$  that is not  in $E(G)$, then replace this edge by a descendant other than $(t-3)(t-2)(t-1)$ with the lowest hierarchy. Repeat this until there is no such an edge.

\end{enumerate}

  Then $G'$ satisfies the following properties:

\begin{enumerate}

\item  The number of edges in $G'$ is the same as the number of edges in $G$.

\item $\lambda(G)=\lambda(G, {\vec x})\le \lambda(G', {\vec x})\le \lambda(G').$

\item $(t-3)(t-2)(t-1) \notin E(G')$.

\item $[t-1]^{(3)}\backslash \{(t-3)(t-2)(t-1)\} \in E(G')$.

\item For any edge in $E(G')$, all its  descendants  other than $(t-3)(t-2)(t-1)$ will be in $E(G')$.

\end{enumerate}

If $G'$ is not left-compressed, then there is an ancestor $uvw$ of $(t-3)(t-2)(t-1)$ such that $uvw\in E(G')$. We claim that $uvw$ must be $(t-3)(t-2)t$. If $uvw$ is not $(t-3)(t-2)t$, then since all  descendants  other than $(t-3)(t-2)(t-1)$ of $uvw$ will be in $E(G')$, then all descendants of $(t-3)(t-1)t$ (other than\\ $(t-3)(t-2)(t-1)$) or
all descendants of $(t-3)(t-2)(t+1)$ (other than $(t-3)(t-2)(t-1)$) will be in $E(G')$. So
all triples in $[t-1]^{(3)}\setminus\{(t-3)(t-2)(t-1)\}$, all triples in the form of $ijt$ (where $ij\in [t-2]^{(2)}$), and all triples in the form of $ij(t+1)$ (where $ij\in [t-2]^{(2)}$) or all triples in the form of  $i(t-1)l$, $1\le i\le t-3$ will be in $E(G')$, then
$$m\ge {t-1 \choose 3}-1 + {t-2 \choose 2}+(t-3)>{t-1 \choose 3} + {t-2 \choose 2},$$ which is a contradiction.
So $uvw$ must be $(t-3)(t-2)t$. Since $m\le {t-1 \choose 3} + {t-2 \choose 2}$ and all the descendants other than
$(t-3)(t-2)(t-1)$ of an edge in $G'$ will be an edge in $G'$,
then there are two possibilities.

Case 1. $E(G')=([t-1]^{(3)}\setminus \{(t-1)(t-2)(t-3)\})\cup\{ijt, ij\in [t-2]^{(2)}\}\cup\{12(t+1)\}.$

Case 2.  $E(G')=([t-1]^{(3)}\setminus \{(t-1)(t-2)(t-3)\})\cup\{ijt, ij\in [t-2]^{(2)}\} .$

Let ${\vec y}:=(y_1, y_2, \ldots, y_n)$ be an optimal weighting of $G'$, where  $n=t+1$ or $n=t$.  We claim that if Case 1 happens, then $y_{t}=y_{t+1}=0$, since $E_{(t-1)t}=E_{t(t+1)}=\emptyset$(by Lemma \ref{LemmaTal5}). If Case 2 happens, then $y_{t}=0$ since $E_{(t-1)t}=\phi$(by Lemma \ref{LemmaTal5}). Hence we can assume that $G$ is left-compressed.
\qed
\bigskip

 \noindent{\em Proof of Lemma \ref{lemma33}} Since $G$ contains the clique of $[t-1]^{(3)}\backslash \{(t-3)(t-2)(t-1)\}$, it is true for $k\leq t$. Next we assume that $k\geq t+1$.

Since $G$ is left-compressed, $1(k-1)k\in E$. Let $b:=\max\{i: i(k-1)k\in E\}$. Since $E$ is left-compressed, then $E_i:=\{1, \ldots, i-1, i+1, \ldots, k\}^{(2)}$, for $1\le i\le b$, and $E_{i\setminus j}=\emptyset$ for $1\le i<j\le b$. Hence, by Remark \ref{r1}(a), we have $x_1=x_2=\cdots=x_b$.

 We define a new feasible weighting ${\vec y}$ for $G$ as follows. Let $y_i=x_i$ for $i\neq k-1, k$, $y_{k-1}=x_{k-1}+x_{k}$ and $y_{k}=0$.

By Lemma \ref{LemmaTal5}(a), $\lambda(E_{k-1}, \vec{x})=\lambda(E_{k}, \vec{x})$, so
\begin{eqnarray}\label{eqa1}
\lambda(G,\vec {y})- \lambda(G,\vec {x})&=&x_{k}(\lambda(E_{k-1}, \vec{x})-x_k\lambda(E_{k(k-1)}, \vec{x}))\nonumber \\
&&-x_{k}(\lambda(E_{k}, \vec{x})-x_{k-1}\lambda(E_{k(k-1)}, \vec{x}))-x_{k-1}x_{k}\lambda(E_{k(k-1)}, \vec{x}) \nonumber \\
&=&x_{k}(\lambda(E_{k-1}, \vec{x})-\lambda(E_{k}, \vec{x}))-x_{k}^2\sum_{i=1}^{b} x_i \nonumber \\
&=&-bx_1x_{k}^2.
\end{eqnarray}
Since $y_{k}=0$ we may remove all edges containing $k$ from $E$ to form a new $3$-graph $\overline{G}:=([k], \overline{E})$ with\\ $\vert \overline{E}\vert:=\vert E\vert-\vert E_{k}\vert$ and $\lambda(\overline{G},\vec {y})=
\lambda(G,\vec {y})$.
We will show that if Lemma \ref{lemma33} fails to hold then there exists a set of edges $F\subset [k-1]^{(3)}\setminus E$ satisfying
\begin{equation}\label{eqa2}
\lambda(F,\vec {y})> bx_1x_{k}^2
\end{equation}
and
\begin{equation}\label{eqa3}
\vert F\vert\le \vert E_{k}\vert.
\end{equation}
Then, using (\ref{eqa1}), (\ref{eqa2}), and (\ref{eqa3}), the $3$-graph $G':=([k], E')$, where $E':=\overline{E}\cup F$, satisfies $\vert E'\vert\le \vert E\vert$ and
\begin{eqnarray*}
\lambda(G',\vec {y})&=&\lambda(\overline{G},\vec {y})+\lambda(F,\vec {y})\\
&>&\lambda(G,\vec {y})+bx_1x_{k}^2 \\
&=&\lambda(G,\vec {x}).
\end{eqnarray*}
Hence $\lambda(G')>\lambda(G)$. Note that  $G'$ still contains $[t-1]^{(3)}\backslash \{(t-3)(t-2)(t-1)\}$ since $G'$ contains all edges in $E\cap [k-1]^{(3)}\supseteq E\cap [t-1]^{(3)}$. If $G'$ does not contains a clique of size $t-1$, note that  $G'$ still contain\\ $[t-1]^{(3)}\backslash \{(t-3)(t-2)(t-1)\}$, it contradicts to $\lambda(G)=\lambda_{(m,t-1)}^{3-}$. If $G'$ contains a clique of size $t-1$, then by Theorem \ref{Corollary 2} $\lambda(G')=\lambda([t-1]^{(3)})$ and consequently $\lambda(G')<\lambda([t-1]^{(3)})$.

We must now construct the set of edges $F$ satisfying (\ref{eqa2}) and (\ref{eqa3}). Applying Remark \ref{r1}(a) by taking $i=1, j=k-1$, we have
$$x_1=x_{k-1}+{\lambda(E_{1\setminus (k-1)}, \vec{x}) \over \lambda(E_{1(k-1)}, \vec{x})}.$$
Let
$C:=[k-2]^{(2)} \setminus E_{k-1}$. Then $\lambda(E_{1\setminus (k-1)}, \vec{x})= x_{k}\sum_{i=b+1}^{k-2}x_i+\lambda(C, \vec{x})$. Applying this and multiplying $bx_{k}^2$ to the above equation (note that $\lambda(E_{1(k-1)}, \vec{x})=\sum_{i=2, i\neq k-1}^{k}x_i$), we have
$$bx_1x_{k}^2=bx_{k-1}x_{k}^2+{bx_{k}^3\sum_{i=b+1}^{k-2}x_i \over \sum_{i=2, i\neq k-1}^{k}x_i}+{bx_{k}^2\lambda(C, \vec{x}) \over \sum_{i=2, i\neq k-1}^{k}x_i}.$$
Since $x_1\ge x_2\ge\cdots\ge x_k$, then
\begin{equation}\label{eqa4}
bx_1x_{k}^2\le bx_{k-1}x_{k}^2(1+{k-(b+2) \over k-3})+{bx_{k}\lambda(C, \vec{x}) \over k-2}.
\end{equation}
Define $\alpha:=\lceil {b\vert C\vert \over k-2} \rceil$ and $\beta:=\lceil b(1+{k-(b+2) \over k-3})\rceil$. Note that $\lceil b(1+{k-(b+2) \over k-3})\rceil \le k-2$ since $b\le k-2$. So $\beta\le k-2$. Let the set $F_1\subset [k-1]^{(3)}\setminus E$ consist of the $\alpha$ heaviest edges in $[k-1]^{(3)}\setminus E$ containing the vertex $k-1$ (note that $\vert [k-2]^{(2)} \setminus E_{k-1}\vert =\vert C\vert \ge \alpha$). Recalling that $y_{k-1}=x_{k-1}+x_k$ we have
$$\lambda(F_1, \vec{y}) \ge {bx_k\lambda(C, \vec{x}) \over k-2}+\alpha x_{k-1}x_{k}^2.$$
So using (\ref{eqa4})
\begin{equation}\label{eqa5}
\lambda(F_1, \vec{y})-bx_{1}x_{k}^2\ge x_{k-1}x_{k}^2(\alpha-\beta).
\end{equation}
We now distinguish two cases.

Case 1. $\alpha>\beta$.

In this case $\lambda(F_1, \vec{y})-bx_{k-1}x_{k}^2>0$ so defining $F:=F_1$ satisfies (\ref{eqa2}). We need to check that $\vert F\vert \le \vert E_k\vert$. Since $E$ is left-compressed, then $[b]^{(2)}\cup \{1, \ldots , b\}\times \{b+1, \ldots, k-1\}\subset E_k$. Hence
\begin{equation}\label{eqa6}
\vert E_k\vert \ge {b[b-1+2(k-1-b)] \over 2} \ge {b(k-1) \over 2}
\end{equation}
since $b\le k-2$.
Recall that $\vert F\vert=\alpha=\lceil {b\vert C\vert \over k-2} \rceil$. Since $C\subset [k-2]^{(2)}$, we have
$\vert C\vert \le {k-2 \choose 2}$. So using (\ref{eqa6}) we obtain
$$\vert F\vert \le \lceil {b(k-3) \over 2} \rceil\le {b(k-1) \over 2} \le \vert E_k\vert.$$
So both (\ref{eqa2}) and (\ref{eqa3}) are satisfied.

Case 2. $\alpha\le\beta$.

Suppose that Lemma \ref{lemma33} fails to hold. So $\vert [k-1]^{(3)} \setminus E\vert\ge k-1\ge \beta +1$ (recall that $\beta\le k-2$). Let $F_2$ consist of any $\beta +1-\alpha$ edges in $[k-1]^{(3)} \setminus (E\cup F_1)$ and define $F:=F_1\cup F_2$. Then since $\lambda(F_2, \vec{y})\ge (\beta+1-\alpha)x_{k-1}^3$ and using (\ref{eqa5}),
$$\lambda(F, \vec{y})-bx_{k-1}x_{k}^2=\lambda(F_1, \vec{y})-bx_{k-1}x_{k}^2+\lambda(F_2, \vec{y})\ge (\beta+1-\alpha)x_{k-1}^3-x_{k-1}x_{k}^2(\beta-\alpha)>0.$$
So (\ref{eqa2}) is satisfied. What remains is to check that $\vert F\vert \le \vert E_k\vert.$ In fact,
$$\vert F\vert =\beta+1\le k-1\le {b(k-1) \over 2} \le \vert E_k\vert$$
when $b\ge 2$. If $b=1$, then,
$$\vert F\vert =\beta+1=3\le k-2={b[b-1+2(k-1-b)] \over 2} \le \vert E_k\vert$$ since $k\ge t\ge 5$.
 \qed

\bigskip

\noindent{\em Proof of Lemma \ref{lemma41}} Let $G$ be a $3$-graph on the vertex set $[n]$ with $m$ edges containing a maximal clique of order $t-2$ such that $\lambda(G)=\lambda_{(m,t-2)}^{3}$. We call such a $G$ an extremal $3$-graph for $m$ and $t-2$. Let $\vec{x}:=(x_{1},x_{2},\cdots,x_{n})$ be an optimal weighting of $G$. We can assume that $x_{i}\geq x_{j}$ when $i<j$ since otherwise we can just relabel the vertices of $G$ and obtain another extremal $3$-graph for $m$ and $t-2$ with an optimal weighting $\vec{x}:=(x_{1},x_{2},\cdots,x_{n})$ satisfying $x_{i}\geq x_{j}$ when $i<j$. Next we obtain a new $3$-graph $G'$ from $G$ by performing the followings

1. If $(t-3)(t-2)(t-1)\in E(G)$, then there is at least one triple in $[t-1]^{(3)}\backslash E(G)$, we replace

$(t-3)(t-2)(t-1)$ by this triple;

2. If an edge in $G$ has a descendant other than $(t-3)(t-2)(t-1)$ that is not in $E(G)$, then replace this edge by a descendant other than $(t-3)(t-2)(t-1)$ with the lowest hierarchy. Repeat this until there is no such an edge.

Then $G'$ satisfies the followings

1.  The number of edges in $G'$ is the same as the number of edges in $G$;

2. $G$ contains the clique $[t-2]^{(3)}$;

3.$\lambda(G)=\lambda(G,\vec{x})\leq \lambda(G',\vec{x})\leq \lambda(G')$;

4.$(t-3)(t-2)(t-1)\notin E(G')$;

5. For any edge in $E(G)$, all its descendants other than $(t-3)(t-2)(t-1)$ will be in $E(G')$.

If $G'$ is not left-compressed, then there is an ancestor $uvw$ of $(t-3)(t-2)(t-1)$ such that $uvw\in G'$ and all the descendant of $uvw$ other than $uvw$ are in $G'$. Hence

$E(G')\supseteq ([t-1]^{(3)}\backslash \{(t-3)(t-2)(t-1)\})\cup \{ijt,ij\in [t-2]^{(2)}\}$.

and $$m\geq {t-1 \choose 3}-1+{t-2 \choose 2}>{t-1 \choose 3}+{t-2 \choose 2}-{t-2 \over 2}.$$
which is a contradiction. Hence $G'$ is left-compressed.
\qed

\bigskip

\noindent{\em  Proof of Lemma \ref{lemma42}} Since $G$ contains the clique of $[t-2]^{(3)}$, it is true for $k\leq t-1$. Assume that $k\geq t$.

Since $G$ is left-compressed, $1(k-1)k\in E$. Let $b:=\max\{i: i(k-1)k\in E\}$. Since $E$ is left-compressed,  $E_i=\{1, \ldots, i-1, i+1, \ldots, k\}^{(2)}$, for $1\le i\le b$, and $E_{i\setminus j}=\emptyset$ for $1\le i<j\le b$. Hence, by Remark \ref{r1}(a), we have $x_1=x_2=\cdots=x_b$.

 We define a new feasible weighting ${\vec y}$ for $G$ as follows. Let $y_i:=x_i$ for $i\neq k-1, k$, $y_{k-1}:=x_{k-1}+x_{k}$ and $y_{k}:=0$.

By Lemma \ref{LemmaTal5}(a), $\lambda(E_{k-1}, \vec{x})=\lambda(E_{k}, \vec{x})$, so
\begin{eqnarray}\label{eqa}
\lambda(G,\vec {y})- \lambda(G,\vec {x})&=&x_{k}(\lambda(E_{k-1}, \vec{x})-x_k\lambda(E_{k(k-1)}, \vec{x}))\nonumber \\
&&-x_{k}(\lambda(E_{k}, \vec{x})-x_{k-1}\lambda(E_{k(k-1)}, \vec{x}))-x_{k-1}x_{k}\lambda(E_{k(k-1)}, \vec{x}) \nonumber \\
&=&x_{k}(\lambda(E_{k-1}, \vec{x})-\lambda(E_{k}, \vec{x}))-x_{k}^2\sum_{i=1}^{b} x_i \nonumber \\
&=&-bx_1x_{k}^2.
\end{eqnarray}
Since $y_{k}=0$ we may remove all edges containing $k$ from $E$ to form a new $3$-graph $\overline{G}:=([k], \overline{E})$ with\\ $\vert \overline{E}\vert:=\vert E\vert-\vert E_{k}\vert$ and $\lambda(\overline{G},\vec {y})=
\lambda(G,\vec {y})$.
We will show that if Lemma \ref{lemma42} fails to hold then there exists a set of edges $F\subset [k-1]^{(3)}\setminus E$ satisfying
\begin{equation}\label{eqb}
\lambda(F,\vec {y})> bx_1x_{k}^2,
\end{equation}
and
\begin{equation}\label{eqc}
\vert F\vert\le \vert E_{k}\vert.
\end{equation}
Then, using (\ref{eqa}), (\ref{eqb}), and (\ref{eqc}), the $3$-graph $G':=([k], E')$, where $E':=\overline{E}\cup F$, satisfies $\vert E'\vert\le \vert E\vert$ and
\begin{eqnarray*}
\lambda(G',\vec {y})&=&\lambda(\overline{G},\vec {y})+\lambda(F,\vec {y})\\
&>&\lambda(G,\vec {y})+bx_1x_{k}^2 \\
&=&\lambda(G,\vec {x}).
\end{eqnarray*}
Hence $\lambda(G')>\lambda(G)$. Note that  $G'$ still contains  the clique $[t-2]^{(3)}$ since $G'$ contains all edges in\\ $E\cap [k-1]^{(3)}\supset [t-2]^{(3)}$. If $G'$ does not contains a clique of size $t-1$, it contradicts to $\lambda(G)=\lambda_{(m,t-2)}^{3}$. If $G'$ contains a clique of size $t-1$, then by Theorem \ref{theorem 1} $\lambda(G')=\lambda([t-1]^{(3)})$ and consequently \\$\lambda(G')<\lambda([t-1]^{(3)})$.

We must now construct the set of edges $F$ satisfying (\ref{eqb}) and (\ref{eqc}). Applying Remark \ref{r1}(a) by taking $i=1, j=k-1$, we have
$$x_1=x_{k-1}+{\lambda(E_{1\setminus (k-1)}, \vec{x}) \over \lambda(E_{1(k-1)}, \vec{x})}.$$
Let
$C:=[k-2]^{(2)} \setminus E_{k-1}$. Then $\lambda(E_{1\setminus (k-1)}, \vec{x})= x_{k}\sum_{i=b+1}^{k-2}x_i+\lambda(C, \vec{x})$. Applying this and multiplying $bx_{k}^2$ to the above equation (note that $\lambda(E_{1(k-1)}, \vec{x})=\sum_{i=2, i\neq k-1}^{k}x_i$), we have
$$bx_1x_{k}^2=bx_{k-1}x_{k}^2+{bx_{k}^3\sum_{i=b+1}^{k-2}x_i \over \sum_{i=2, i\neq k-1}^{k}x_i}+{bx_{k}^2\lambda(C, \vec{x}) \over \sum_{i=2, i\neq k-1}^{k}x_i}.$$
Since $x_1\ge x_2\ge\cdots\ge x_k$, then
\begin{equation}\label{eqd}
bx_1x_{k}^2\le bx_{k-1}x_{k}^2(1+{k-(b+2) \over k-3})+{bx_{k}\lambda(C, \vec{x}) \over k-2}.
\end{equation}
Define $\alpha:=\lceil {b\vert C\vert \over k-2} \rceil$ and $\beta:=\lceil b(1+{k-(b+2) \over k-3})\rceil$. Note that $\lceil b(1+{k-(b+2) \over k-3})\rceil \le k-2$ since $b\le k-2$. So $\beta\le k-2$. Let the set $F_1\subset [k-1]^{(3)}\setminus E$ consist of the $\alpha$ heaviest edges in $[k-1]^{(3)}\setminus E$ containing the vertex $k-1$ (note that $\vert [k-2]^{(2)} \setminus E_{k-1}\vert =\vert C\vert \ge \alpha$). Recalling that $y_{k-1}=x_{k-1}+x_k$ we have
$$\lambda(F_1, \vec{y}) \ge {bx_k\lambda(C, \vec{x}) \over k-2}+\alpha x_{k-1}x_{k}^2.$$
So using (\ref{eqd})
\begin{equation}\label{eqe}
\lambda(F_1, \vec{y})-bx_{1}x_{k}^2\ge x_{k-1}x_{k}^2(\alpha-\beta).
\end{equation}
We now distinguish two cases.

Case 1. $\alpha>\beta$.

In this case $\lambda(F_1, \vec{y})-bx_{k-1}x_{k}^2>0$ so defining $F:=F_1$ satisfies (\ref{eqb}). We need to check that $\vert F\vert \le \vert E_k\vert$. Since $E$ is left-compressed, then $[b]^{(2)}\cup \{1, \ldots , b\}\times \{b+1, \ldots, k-1\}\subset E_k$. Hence
\begin{equation}\label{eqf}
\vert E_k\vert \ge {b[b-1+2(k-1-b)] \over 2} \ge {b(k-1) \over 2}
\end{equation}
since $b\le k-2$.
Recall that $\vert F\vert=\alpha=\lceil {b\vert C\vert \over k-2} \rceil$. Since $C\subset [k-2]^{(2)}$, we have
$\vert C\vert \le {k-2 \choose 2}$. So using (\ref{eqd}) we obtain
$$\vert F\vert \le \lceil {b(k-3) \over 2} \rceil\le {b(k-1) \over 2} \le \vert E_k\vert.$$
So both (\ref{eqb}) and (\ref{eqc}) are satisfied.

Case 2. $\alpha\le\beta$.

Suppose that Lemma \ref{lemma42} fails to hold. So $\vert [k-1]^{(3)} \setminus E\vert\ge k-1\ge \beta +1$ (recall that $\beta\le k-2$). Let $F_2$ consist of any $\beta +1-\alpha$ edges in $[k-1]^{(3)} \setminus (E\cup F_1)$ and define $F:=F_1\cup F_2$. Then since $\lambda(F_2, \vec{y})\ge (\beta+1-\alpha)x_{k-1}^3$ and using (\ref{eqe}),
$$\lambda(F, \vec{y})-bx_{k-1}x_{k}^2=\lambda(F_1, \vec{y})-bx_{k-1}x_{k}^2+\lambda(F_2, \vec{y})\ge (\beta+1-\alpha)x_{k-1}^3-x_{k-1}x_{k}^2(\beta-\alpha)>0.$$
So (\ref{eqb}) is satisfied. What remains is to check that $\vert F\vert \le \vert E_k\vert.$ In fact,
$$\vert F\vert =\beta+1\le k-1\le {b(k-1) \over 2} \le \vert E_k\vert$$
when $b\ge 2$. If $b=1$, then applying (\ref{eqe}),
$$\vert F\vert =\beta+1=3\le k-2={b[b-1+2(k-1-b)] \over 2} \le \vert E_k\vert$$ since $k\ge t\ge 5$.
 \qed

\bigskip

\noindent{\em  Proof of Lemma \ref{lemma43}} Let $b:=\max\{i: i(t-1)t\in E\}$. Since $E$ is left-compressed, then

$E_i=\{1, \ldots, i-1, i+1, \ldots, t\}^{(2)}, $ for $1\le i\le b$ and $E_{i\setminus j}=\emptyset$ for $1\le i<j\le b$.

Hence, by Remark \ref{r1}(a), we have $x_1=x_2=\cdots=x_b$. Consider a new weighting for $G$, ${\vec z}:=(z_1, z_2, \ldots, z_t)$ given by $z_i:=x_i$ for $i\neq t-1, t$, $z_{t-1}:=0$ and $z_t:=x_{t-1}+x_t$. By Lemma \ref{LemmaTal5}(a), $\lambda(E_{t-1}, \vec{x})=\lambda(E_{t}, \vec{x})$, so
\begin{eqnarray}\label{eq10}
\lambda(G,\vec {z})- \lambda(G,\vec {x})&=&x_{t-1}(\lambda(E_{t}, \vec{x})-\lambda(E_{t-1}, \vec{x}))-x_{t-1}^2\sum_{i=1}^b x_i
=-bx_1x_{t-1}^2.
\end{eqnarray}

Since $z_{t-1}=0$ we may remove all edges containing $t-1$ from $E$ to form a new $3$-graph $\overline{G}:=([t], \overline{E})$ with $\vert \overline{E}\vert:=\vert E\vert-\vert E_{t-1}\vert$ and $\lambda(\overline{G},\vec {z})=
\lambda(G,\vec {z})$.

If $|[t-2]^{(2)} \backslash E_t| > b$, we will show that there exists a set of edges $F\subset \{1, ..., t-2, t\}^{(3)}\setminus E$ satisfying
\begin{equation}\label{eq11}
\lambda(F,\vec {z})> bx_1x_{t-1}^2.
\end{equation}

Then using (\ref{eq10}) and (\ref{eq11}), the $3$-graph $G':=([t], E')$, where $E':=\overline{E}\cup F$, satisfies $\lambda(G', \vec {z}))>\lambda(G)$. Since $\vec {z}$ has only $t-1$ positive weights, then $\lambda(G', \vec {z}))\le \lambda([t-1]^{(3)})$, and consequently
$$\lambda(G)<\lambda([t-1]^{(3)}).$$
We must now construct the set of edges $F$. Since $G$ is left-compressed, applying  Remark \ref{r1}(a) by taking $i=1$, $j=t$, we get
$$x_1=x_t+{\lambda(E_{1 \setminus t}, \vec{x}) \over \lambda(E_{1t}, \vec{x})}.$$
Let
$D:=[t-2]^{(2)} \setminus E_t$. Then $\lambda(E_{1 \setminus t}, \vec{x})= x_{t-1}\sum_{i=b+1}^{t-2}x_i+\lambda(D, \vec{x})$. Applying this and  multiplying $bx_{t-1}^2$ to the above equation (note that $\lambda(E_{1t}, \vec{x})=\sum_{i=2}^{t-1}x_i$), we have
$$bx_1x_{t-1}^2=bx_tx_{t-1}^2+{bx_{t-1}^3\sum_{i=b+1}^{t-2}x_i \over \sum_{i=2}^{t-1}x_i}+{bx_{t-1}^2\lambda(D, \vec{x}) \over \sum_{i=2}^{t-1}x_i}.$$
Let $c:={\sum_{i=b+1}^{t-2}x_i \over \sum_{i=2}^{t-1}x_i}$ and  $d:={bx_{t-1} \over \sum_{i=2}^{t-1}x_i}$. Then
\begin{equation}\label{eq13}
bx_1x_{t-1}^2=bx_tx_{t-1}^2+bcx_{t-1}^3+dx_{t-1}\lambda(D, \vec{x}).
\end{equation}
Let $F$ consist of those edges in $\{1, ..., t-2, t\}^{(3)}\setminus E$ containing the vertex $t$. Then
\begin{equation}\label{eq132}
\lambda(F,\vec {z})=(x_{t-1}+x_l)\lambda(D, \vec{x}).
\end{equation}
Since $|[t-2]^{(2)} \backslash E_t| > b$, then
\begin{equation}\label{eq133}
\lambda(D, \vec{x}) > bx_{t-1}^2.
\end{equation}
Applying equations (\ref{eq13}), (\ref{eq132}), and (\ref{eq133}), we  get
\begin{eqnarray*}
\lambda(F,\vec {z})-bx_1x_{t-1}^2&=&(x_{t-1}+x_{t})\lambda(D, \vec{x})-bx_{t}x_{t-1}^2-bcx_{t-1}^3-dx_{t-1}\lambda(D, \vec{x})\\
&=&[(1-d)x_{t-1}+x_{t}]\lambda(D, \vec{x})-bx_{t}x_{t-1}^2-bcx_{t-1}^3\\
&>&[(1-d)x_{t-1}+x_{t}]bx_{t-1}^2-bx_{t}x_{t-1}^2-bcx_{t-1}^3\\
&=& bx_{t-1}^3(1-d-c)\ge 0.
\end{eqnarray*}
since $$c+d={\sum_{i=b+1}^{t-2}x_i+bx_{t-1} \over \sum_{i=2}^{t-1}x_i}\le 1. $$

Let $G':=([t], \overline{E} \cup F)$, then
$\lambda(G', \vec{z})=\lambda(G,\vec{z}))+\lambda(F,\vec{z})=\lambda(G,\vec{x})-bx_1x_{t-1}^2+\lambda(F,\vec{z}) > \lambda(G,\vec{x})$. On the other hand, since $\vec{z}$ has only $t-1$ positive weights, then $\lambda(G', \vec{z}) < \lambda([t-1]^{(3)})$. \qed

\bigskip

\noindent{\em Proof of Lemma \ref{lefcom}} Let $m$ and $t$ be positive integers satisfying ${t-1 \choose r} \le m \le {t \choose r}-1.$ Let $G:=(V,E)$ be an $r$-graph on vertex set $V:=[n]$ with $m$ edges containing a clique  of size $t-2$ such that
 $\lambda(G)=\lambda_{(m,t-2)}^{r}$. We call such a $G$ an extremal $r$-graph for $m$ and $t-2$. Let  $\vec{x}:=(x_{1},x_{2},\ldots ,x_{n})$ be an optimal weighting of $G$. We can assume that $x_{i}\geq x_{j}$ when $i<j$ since otherwise we can just relabel the vertices of $G$ and obtain another extremal $r$-graph for $m$ and $t-2$ with an optimal $\vec{x}:=(x_{1},x_{2},\ldots ,x_{n})$ satisfying $x_{i}\geq x_{j}$ when $i<j$. If $G$ is not left-compressed, then there is an edge whose ancestor is not an edge. Replace all those edges by its available ancestor with the highest hierarchy, then we get a left-compressed $r$-graph $G'$ which contains the clique $[t-2]^{(r)}$ and $\lambda(G',\vec{x})\geq \lambda(G,\vec{x})$.\qed

\bigskip

\noindent{\em Proof of Lemma \ref{Lemma51}}  Let $b:=\max\{i: i(t-1)t\in E\}$. Since $E$ is left-compressed, then

$E_i=\{1, \ldots, i-1, i+1, \ldots, t\}^{(2)}$, for $1\le i\le b$, and $E_{i\setminus j}=\emptyset$ for $1\le i<j\le b$.

Hence, by Remark \ref{r1}(a), we have $x_1=x_2=\cdots=x_b$. We define a new feasible weighting ${\vec y}$ for $G$ as follows. Let $y_i:=x_i$ for $i\neq t-1, t$, $y_{t-1}:=x_{t-1}+x_{t}$ and $y_{t}:=0$.

By Lemma \ref{LemmaTal5}(a), $\lambda(E_{t-1}, \vec{x})=\lambda(E_{t}, \vec{x})$, so
\begin{eqnarray}\label{eqb1}
\lambda(G,\vec {y})- \lambda(G,\vec {x})&=&x_{t}(\lambda(E_{t-1}, \vec{x})-x_t\lambda(E_{t(t-1)}, \vec{x}))\nonumber \\
&&-x_{t}(\lambda(E_{t}, \vec{x})-x_{t-1}\lambda(E_{(t-1)t}, \vec{x}))-x_{t-1}x_{t}\lambda(E_{(t-1)t}, \vec{x}) \nonumber \\
&=&x_{t}(\lambda(E_{t-1}, \vec{x})-\lambda(E_{t}, \vec{x}))-x_{t}^2\sum_{i=1}^{b} x_i \nonumber \\
&=&-bx_1x_{t}^2.
\end{eqnarray}
Since $y_{t}=0$ we may remove all edges containing $t$ from $E$ to form a new $3$-graph

 $\overline{G}:=([t], \overline{E})$ with $\vert \overline{E}\vert:=\vert E\vert-\vert E_{t}\vert$ and $\lambda(\overline{G},\vec {y})=\lambda(G,\vec {y})$.

We will show that if $|[t-1]^{(3)}\setminus E|\geq t-2$ then there exists a set of edges $F\subset [t-1]^{(3)}\setminus E$ satisfying
\begin{equation}\label{eqb2}
\lambda(F,\vec {y})\geq bx_1x_{t}^2,
\end{equation}
Then, using (\ref{eqb1}), (\ref{eqb2}), the $3$-graph $G':=([t], E')$, where $E':=\overline{E}\cup F$, satisfies
\begin{eqnarray*}
\lambda(G',\vec {y})&=&\lambda(\overline{G},\vec {y})+\lambda(F,\vec {y})\\
&\geq &\lambda(G,\vec {y})+bx_1x_{t}^2 \\
&=&\lambda(G,\vec {x}).
\end{eqnarray*}
Since $\vec{y}$ has only $t-1$ positive weights, then $ \lambda(G')\leq \lambda([t-1]^{(3)})$, and consequently $ \lambda(G)\leq \lambda([t-1]^{(3)})$.

We must now construct the set of edges $F$ satisfying (\ref{eqb2}). Applying Remark \ref{r1}(a) by taking $i=1, j=t-1$, we have
$$x_1=x_{t-1}+{\lambda(E_{1\setminus (t-1)}, \vec{x}) \over \lambda(E_{1(t-1)}, \vec{x})}.$$
Let
$C:=[t-2]^{(2)} \setminus E_{t-1}$. Then $\lambda(E_{1\setminus (t-1)}, \vec{x})= x_{t}\sum_{i=b+1}^{t-2}x_i+\lambda(C, \vec{x})$. Applying this and multiplying $bx_{t}^2$ to the above equation (note that $\lambda(E_{1(t-1)}, \vec{x})=\sum_{i=2, i\neq t-1}^{t}x_i$), we have
$$bx_1x_{t}^2=bx_{t-1}x_{t}^2+{bx_{t}^3\sum_{i=b+1}^{t-2}x_i \over \sum_{i=2, i\neq t-1}^{t}x_i}+{bx_{t}^2\lambda(C, \vec{x}) \over \sum_{i=2, i\neq t-1}^{t}x_i}.$$
Since $x_1\ge x_2\ge\cdots\ge x_t$, then
\begin{equation}\label{eqb3}
bx_1x_{t}^2\le bx_{t-1}x_{t}^2(1+{t-(b+2) \over t-3})+{bx_{t}\lambda(C, \vec{x}) \over t-2}.
\end{equation}
Define $\alpha:=\lceil {b\vert C\vert \over t-2} \rceil$ and $\beta:=\lceil b(1+{t-(b+2) \over t-3})\rceil$. Note that since
$b\le t-2$. So $\beta\le t-2$. Let the set \\$F_1\subset [t-1]^{(3)}\setminus E$ consist of the $\alpha$ heaviest edges in $[t-1]^{(3)}\setminus E$ containing the vertex $t-1$ (note that \\$\vert [t-2]^{(2)} \setminus E_{t-1}\vert =\vert C\vert \ge \alpha$). Recalling that $y_{t-1}=x_{t-1}+x_t$ we have
$$\lambda(F_1, \vec{y}) \ge {bx_t\lambda(C, \vec{x}) \over t-2}+\alpha x_{t-1}x_{t}^2.$$
So using (\ref{eqb3})
\begin{equation}\label{eqb4}
\lambda(F_1, \vec{y})-bx_{1}x_{t}^2\ge x_{t-1}x_{t}^2(\alpha-\beta).
\end{equation}

If $\alpha>\beta$, then $\lambda(F_1, \vec{y})-bx_{t-1}x_{t}^2>0$. So defining $F:=F_1$ satisfies (\ref{eqb2}).

Assume $\alpha\le\beta$. Suppose that $|[t-1]^{(3)}\backslash E|\geq t-2$. So $\vert [t-1]^{(3)} \setminus E\vert\ge t-2\ge \beta$ (recall that $\beta\le t-2$). Let $F_2$ consist of any $\beta-\alpha$ edges in $[t-1]^{(3)} \setminus (E\cup F_1)$ and define $F:=F_1\cup F_2$. Then since $\lambda(F_2, \vec{y})\ge (\beta-\alpha)x_{t-1}^3$ and using (\ref{eqb3})
$$\lambda(F, \vec{y})-bx_{t-1}x_{t}^2=\lambda(F_1, \vec{y})-bx_{t-1}x_{t}^2+\lambda(F_2, \vec{y})\ge (\beta-\alpha)x_{t-1}^3-x_{t-1}x_{t}^2(\beta-\alpha)\geq 0.$$
This proves Lemma \ref{Lemma51}. \qed

\section{Conclusions}
At this moment, we are not able to extend the arguments in this paper to verify Conjectures \ref{conjecture1}, \ref{conjecture2}, and \ref{conjecture} for more general cases. When $r\ge 4$,  the computation is more complex.  If  there is some technique to overcome this difficulty, then the idea used in proving Theorem \ref{theorem31}   can  be used to   improve our results much further.
\bigskip\\
{\bf Acknowledgments} We thank two anonymous referees and the editor for  helpful and insightful  comments. This research is partially  supported by National Natural Science Foundation of China (No. 11271116).

\bibliographystyle{spmpsci}      
\bibliographystyle{unsrt}


\end{document}